\documentclass[12pt]{article}
\input{epsf}

\newtheorem{theorem}{Theorem}
\newtheorem{lemma}{Lemma}

\newenvironment{proof}{{\bf Proof.}}{\hfill\rule{2mm}{2mm}}

\newtheorem{remarka}{Remark}
\newenvironment{remark}{\begin{remarka}\rm}{\end{remarka}}
\newtheorem{prelem}{{\bf Theorem}}

\newenvironment{lem}{\begin{prelem}{\hspace{-0.5
               em}{\bf.}}}{\end{prelem}}

\newtheorem{preconj}{{\bf Conjecture}}

\newenvironment{conj}{\begin{preconj}{\hspace{-0.5
               em}{\bf.}}}{\end{preconj}}

\def\PrConst {180}
\def\MAXDEG {290}
\def\MINDEG {20}
\def\UDEG {100}

\def\UC {{\mathit{UC}}}
\def\Prob {{\rm \bf Pr}}
\def\Ex {{\rm \bf E}}
\def\deg {{\rm deg}}
\def\avd {{\em avd}}

\def\final{{\mathit final}}
\def\BIN {{\mathit BIN}}

\title{\bf $\Delta+300$ is a Bound on the Adjacent Vertex Distinguishing Edge Chromatic Number}
\author{ Hamed Hatami \\
{\small\it Department of Computer Science}\\
{\small University of Toronto} \\
{\small e-mail: hamed@cs.toronto.edu}}
\date{}
\begin{document}
\maketitle

\begin{abstract}
An adjacent vertex distinguishing edge-coloring or an
\avd-coloring of a simple graph $G$ is a proper edge-coloring of
$G$ such that no pair of adjacent vertices meets the same set of
colors. We prove that every graph with maximum degree $\Delta$ and
with no isolated edges has an \avd-coloring with at most
$\Delta+300$ colors, provided that $\Delta >10^{20}$.
\end{abstract}

\noindent {{\sc AMS Subject Classification:} \quad  05C15}
\newline
{{\sc Keywords:} adjacent vertex distinguishing edge coloring,
vertex distinguishing edge coloring, adjacent strong edge
coloring, strong edge coloring.

\section{Introduction}
We follow~\cite{west} for terminologies and notations not defined
here. For every vertex $v$ in a graph $G$, $\deg_G(v)$ or
$\deg(v)$ when there is no ambiguity denotes the degree of $v$ in
the graph $G$. For every partial edge coloring $c$ of a graph $G$,
and every vertex $v$ of $G$, let $S_c(v)$ denote the set of the
colors incident to $v$.

In a partial edge coloring $c$ a vertex $u$ is called {\sf
distinguishable}, if for every vertex $v$ adjacent to $u$, $S_c(u)
\neq S_c(v)$. An edge coloring $c$ is called {\sf adjacent vertex
distinguishing} or an {\sf \avd-coloring} if every vertex is
distinguishable. A {\sf $k$-\avd-coloring} is an \avd-coloring
using at most $k$ colors. It is clear that every graph with
isolated edges does not have any \avd-coloring. The {\sf
\avd-chromatic number} of a graph $G$, the minimum number of
colors in an \avd-coloring of $G$, is defined for every graph $G$
without any isolated edge. Adjacent vertex distinguishing edge
colorings are studied in~\cite{Akbari,Balister2,ghandehari,Zhang},
where different names such as {\sf adjacent strong edge
coloring}~\cite{Zhang} and {\sf $1$-strong edge
coloring}~\cite{Akbari, ghandehari} are used to refer to an
\avd-coloring. Adjacent vertex distinguishing edge colorings are
related to vertex distinguishing edge colorings in which the
condition $S_c(u) \neq S_c(v)$ holds for every pair of vertices
$u$ and $v$, not necessarily adjacent. This concept has been
studied in many papers, (see for
example~\cite{Balister,Burris,favaron,Hornak2}).

Another interesting problem arises when we drop the condition that
the edge coloring is proper and allow the incident edges to have
the same colors. The following theorem is proved in~\cite{Luczak}
by Karo\'{n}ski, {\L}uczak, and Thomason.

\begin{lem}
{\bf \cite{Luczak}} \label{weight} There exists a finite set of
real numbers which can be used to weight the edges of any graph
with no isolated edges so that adjacent vertices have different
sums of incident edge weights.
\end{lem}
It is easy to observe that Theorem~\ref{weight} is equivalent to
the following theorem.

\begin{lem} {\rm [}restatement of Theorem~\ref{weight}{\rm]}
There is a finite set which can be used to color the edges of any
graph with no isolated edges so that adjacent vertices meet
different multisets (i.e. duplicate elements are allowed) of
colors.
\end{lem}

This shows that by dropping the condition of being proper from the
definition of \avd-coloring, a constant number of colors would be
sufficient. Obviously, when the edge colorings are required to be
proper this is not the case. The following conjecture was made
in~\cite{Zhang}.

\begin{conj}
{\bf \cite{Zhang}} The \avd-chromatic number of every simple
connected graph $G$ such that $G \neq C_5$ (the cycle of size $5$)
and $G\neq K_2$ is at most $\Delta+2$.
\end{conj}

In~\cite{Balister2} this conjecture is verified for bipartite
graphs and graphs of maximum degree $3$, and the following bound
has been proved for general graphs.

\begin{lem}
{\bf \cite{Balister2}} If $G$ is a graph with no isolated edges,
then the \avd-chromatic number of $G$ is at most $\Delta+O(\log
\chi(G))$.
\end{lem}

This was the best known bound so far. We asymptotically improve it
by proving the following theorem.

\begin{theorem}
\label{main} If $G$ is a graph with no isolated edges and maximum
degree $\Delta>10^{20}$, then the \avd-chromatic number of $G$ is
at most $\Delta+300$.
\end{theorem}

We will use the probabilistic method to prove Theorem~\ref{main}.
The following tools of the probabilistic method will be used
several times (see~\cite{molloy}).

\begin{lem} {\rm [}The Symmetric Lov\'{a}sz Local Lemma{\rm]} Let $A_1, \ldots,
A_n$ be events in a probability space, let for all $i$,
$\Prob[A_i]=p$ and the event $A_i$ is mutually independent of all
events but at most $d$ events. If $4pd<1$, then $\Prob[\wedge \bar
A_i]>0$.
\end{lem}

\begin{lem} {\rm [}The Chernoff Bound{\rm]} Suppose $\BIN(n,p)$ is
the sum of $n$ independent Bernoulli variables each occurring with
probability $p$. Then for any $0<t\le np$
$$\Prob[|\BIN(n,p)-np|>t]<2e^{\frac{-t^2}{3np}}.$$
\end{lem}

\begin{lem} {\rm [}Talagrand's Inequality{\rm]}
Let $X$ be a non-negative random variable, not identically $0$,
which is determined by $n$ independent trials $T_1,\ldots,T_n$,
and satisfying the following for some $c,r>0$

\begin{enumerate}
\item Changing the outcome of any one trial can affect $X$ by at
most $c$, and

\item For any $s$ and any outcome of trials, if $X \ge s$, then
there is a set of at most $rs$ trials whose outcomes certify that
$X \ge s$,
\end{enumerate}
then for any $0 \le t \le \Ex[x]$,

$$\Prob[|X-\Ex[x]|>t + 60 c\sqrt{r\Ex[x]}] \le 4
e^{-\frac{t^2}{8c^2r\Ex[X]}}.$$
\end{lem}
\begin{remark}
Suppose that we have the upperbound $\Ex[X] \le k$, and the
conditions of Talagrand's Inequality hold for $X$. Then we might
apply Talagrand's Inequality to the random variable
$Y=X+k-\Ex[X]$, and obtain the inequality
$$\Prob[X>k+t + 60 c\sqrt{rk}] \le \Prob[Y>k+t + 60 c\sqrt{rk}] \le 4
e^{-\frac{t^2}{8c^2rk}}.$$

Similar arguments hold for the Chernoff bound, and also for the
cases that we know a lower bound for $\Ex[X]$ and want to apply
Talagrand's Inequality or the Chernoff bound.
\end{remark}
We will also use the following well known inequalities.
\begin{equation}
\label{ineq1}(\frac{a}{b})^b \le {a \choose b} = {a \choose a-b}
\le (\frac{ea}{b})^b,
\end{equation}
where $a$ and $b$ are natural numbers.

In Section~2 we present the proof of Theorem~\ref{main}.

\section{Proof of Theorem~\ref{main}}
The proof is probabilistic, and consists of three steps. The first
step is similarly used in~\cite{ghandehari}. Let $G$ be a graph
with no isolated edges.
\begin{itemize}
\item In the first step we construct a loopless multigraph $G'$
with multiplicity at most $2$ so that $\Delta(G')
=\Delta(G)=\Delta$, and if there exists a $p$-\avd-coloring for
$G'$, then $G$ is also $p$-\avd-colorable.

\item Let $H$ be the subgraph of $G'$ induced by all vertices $v$
such that $\deg(v)<\frac{\Delta}{3}$. In the second step we give a
$(\Delta+300)$-edge coloring $c$ of $G'$ such that for every two
adjacent vertices $u,v \in V(G')-V(H)$, $S_c(u) \neq S_c(v)$.

\item In the third step we modify $c$ by recoloring some edges in
$E(H)$, and obtain a $(\Delta+300)$-\avd-coloring, $c_{\final}$ of
$G'$.
\end{itemize}

\subsection{First Step}
For every edge $e=uv$ in $G$ such that $\deg(u),\deg(v) <
\frac{\Delta}{3}$ and none of $u$ and $v$ has any other neighbor
of degree less than $\frac{\Delta}{3}$, contract the edge $e$,
i.e. remove $e$ and identify $u$ and $v$ to a single vertex. Call
this new multigraph $G'$. Obviously, $\Delta(G')=\Delta(G)$, and
the multiplicity of $G'$ is at most $2$. Note that $H$, the
subgraph of $G'$ induced by all vertices $v$ such that
$\deg(v)<\frac{\Delta}{3}$, does not have any isolated edges.

Suppose that $G'$ has a $p$-\avd-coloring $\sigma'$. Every edge of
$G'$ has a corresponding edge in $G$. For every edge $e$ of $G'$,
color its corresponding edge in $G$ by $\sigma'(e)$. Now all edges
of $G$ are colored except those edges which were contracted in the
process of obtaining $G'$ from $G$. For every such edge $uv$ in
$G$, since $\deg(v)+\deg(u) < \Delta$, there exist some available
colors for $uv$. Color $uv$ with one of these available colors.
Trivially $u$ and  $v$ meet different sets of colors and none of
them has any other neighbor of the same degree. This implies that
the obtained edge coloring of $G$ is in fact an \avd-coloring.

\subsection{Second Step}
An {\sf unused edge} in a partial edge coloring $c$ of a graph $G$
is an edge which is not colored by $c$. We will refer to the graph
$U_{c}$ induced by all unused edges as the {\sf unused graph} of
$c$. The degree of a vertex in $U_c$ is called its {\sf unused
degree}. We say that a color $x$ is {\sf available} for an unused
edge $uv$, if $x$ does not appear on any incident edge to $uv$.

As it is mentioned before, $H$ is the subgraph of $G'$ induced by
all vertices $v$ such that $\deg(v)<\frac{\Delta}{3}$. It is
trivial that $H$ has no multiple edges. We start with an arbitrary
$\Delta+2$ edge coloring $\sigma$ of $G'$ (which is guaranteed to
exist by Vizing's Theorem), and then we apply a two-phased
procedure to modify $\sigma$, and obtain an edge coloring $c$ such
that for every two adjacent vertices $u,v \in V(G')-V(H)$, $S_c(u)
\neq S_c(v)$. The first phase is the following.
\\
\\
\noindent {\bf Phase I}
\begin{enumerate}
\item Uncolor each edge $e \in E(G')-E(H)$ with probability
$\frac{\PrConst}{\Delta}$.

\item For every vertex $v$ which has unused degree greater than
$\MAXDEG$, recover the color of all unused edges which are
incident to it. We say that $v$ is {\sf recovered}.
\end{enumerate}

Let $\sigma_1$ be the partial coloring obtained after applying
Phase~{\bf I}. We define the following sets.

\begin{itemize}
\item $\UC_v$ (uncolored) is the set of the edges which are
incident to $v$ and uncolored in Phase~{\bf I}~(1). Note that
these edges are not necessarily unused in $\sigma_1$ because their
colors might be recovered in Phase~{\bf I}~(2).

\item $R$ is the set of all recovered vertices.

\item $Q$ is the set of all vertices $v$ such that $|\UC_v|
<\MINDEG$.

\item $T$ is the set of all vertices $v \in V(G')-V(H)$ such that
there exists an edge $vw \in \UC_v$ where $w$ is recovered.

\item $L$ is the set of all vertices $v \in V(G')-V(H)$ such that
$\deg_{U_{\sigma_1}}(v) < \MINDEG$, where $U_{\sigma_1}$ is the
unused graph of $\sigma_1$.
\end{itemize}

Since we use randomness in Phase~{\bf I}, intuitively after
applying this phase if a vertex has a large unused degree, then
with a high probability it is distinguishable. However there are
vertices in $V(G')-V(H)$ which have small unused degrees which we
denote them by $L$. Later in Lemma~\ref{phaseI} we will show that
with a positive probability the vertices in $L$ are rare and
well-distributed in the graph. In fact we will prove that with
positive probability:
\begin{enumerate}
\item[$(a)$:] For every vertex $v \in V(G')-V(H)$, we have $|N(v)
\cap L| \le \frac{\Delta}{\UDEG}$,

\item[$(b)$:] For every two adjacent vertices $u,v \in V(G')-V(H)$
such that $\deg(u)=\deg(v)$ and $u \not\in L$, we have
$|S_{\sigma_1}(u) \bigtriangleup S_{\sigma_1}(v)| \ge 10$, where
for every two sets $A$ and $B$, $A \bigtriangleup B=(A\setminus B)
\cup (B \setminus A)$.
\end{enumerate}

Note that $L$ is a subset of the union of $R$, $Q$, and $T$. The
following three lemmas provide the technical details needed to
prove Lemma~\ref{phaseI}.

\begin{lemma}
\label{probs} For every vertex $v \in V(G')-V(H)$, we have
\begin{enumerate}
\item[$(a):$] $\Prob[v \in R] \le \frac{1}{1000}$.

\item[$(b):$] $\Prob[v \in Q] \le \frac{1}{1000}$.

\item[$(c):$] $\Prob[v \in T] \le \frac{1}{1000}$.
\end{enumerate}
\end{lemma}
\begin{proof}

\noindent (a) Since $|\UC_v|$ has Bernoulli distribution, knowing
that $\frac{\Delta}{3} \le \deg(v) \le \Delta$, the Chernoff bound
implies that
$$\Prob[v \in R]=\Prob[|\UC_v|>\MAXDEG]
\le 2e^{-\frac{(\MAXDEG-\PrConst)^2}{3\times \PrConst}} \le
\frac{1}{1000}.$$
\noindent (b) The worst case is when $\deg(v)=\frac{\Delta}{3}$.
Then the Chernoff bound implies that
$$\Prob[v \in Q]=\Prob[|\UC_v|<\MINDEG]
\le 2e^{-\frac{(\PrConst/3-\MINDEG)^2}{\PrConst}} \le
\frac{1}{1000}.$$
\noindent (c) Consider an edge $vw$, where $v \in V(G')-V(H)$ and
$w \in V(G')$. First we prove an upper bound for $\Prob[vw \in
\UC_v \wedge w \in R]$. The Chernoff bound implies that
$$\Prob[w \in R|vw \in \UC_v] \le 2e^{-\frac{(\MAXDEG-1-\PrConst)^2}{3\times \PrConst}} \le
\frac{1}{10^6}.$$
Hence
$$\Prob[vw \in \UC_v \wedge w \in R] \le \frac{1}{10^6}
\frac{\PrConst}{\Delta} \le \frac{1}{10^3\Delta}.$$
Since $v$ has at most $\Delta$ neighbors, we have
$$\Prob[v \in T] \le  \frac{1}{10^3\Delta} \times \Delta \le \frac{1}{1000}.$$
\end{proof}

\begin{lemma}
\label{intersectL} For every vertex $v \in V(G')-V(H)$,
$$\Prob[|N(v) \cap L|
> \Delta/100] \le \frac{1}{\Delta^7}.$$
\end{lemma}
\begin{proof}
 First notice that
$$|N(v) \cap L| \le |N(v) \cap R| + |N(v) \cap T|+|N(v) \cap Q|.$$
Lemma~\ref{probs} implies that
\begin{itemize}
\item $\Ex\left[|N(v) \cap R| \right] \le \frac{\Delta}{1000}$,

\item $\Ex\left[|N(v) \cap Q| \right] \le \frac{\Delta}{1000}$,
and

\item $\Ex\left[|N(v) \cap T| \right] \le \frac{\Delta}{1000}$.
\end{itemize}

For edges $uv \in E(G')-E(H)$, consider the independent Bernoulli
trials $T_{uv}$  where the outcome of $T_{uv}$ determines whether
$uv$ is uncolored in Phase~{\bf I}~(1) or not. Next we apply
Talagrand's Inequality to prove the following three claims.

\begin{itemize}
\item Claim 1: $\Prob[|N(v) \cap R|
> \Delta/300] \le \frac{1}{3\Delta^{7}}$.

\item Claim 2: $\Prob[|N(v) \cap T|
> \Delta/300] \le \frac{1}{3\Delta^{7}}$.

\item Claim 3: $\Prob[|N(v) \cap Q|
> \Delta/300] \le \frac{1}{3\Delta^{7}}$.
\end{itemize}
To see that these three claims imply Lemma~\ref{intersectL}, we
observe that if $|N(v) \cap L| > \frac{\Delta}{100}$, then either
$|N(v) \cap R|>\frac{\Delta}{300}$, $|N(v) \cap T|
>\frac{\Delta}{300}$, or $|N(v) \cap Q|>\frac{\Delta}{300}$.

\noindent {\bf Proof of Claims 1 and 2:} Changing the outcome of
each trial affects $|N(v) \cap R|$ by at most $2$, and every
assignment to trials that results $|N(v) \cap R| \ge k$ can be
certified by the outcome of $\MAXDEG k$ trials. Hence Talagrand's
Inequality implies that
$$\Prob[|N(v) \cap R| >\frac{\Delta}{1000}+ t+ 120 \sqrt{\MAXDEG
\frac{\Delta}{1000}}]\le 4 e^{- 1000t^2 / 8 \times 4 \times
\MAXDEG \Delta}.$$
Since changing the outcome of each trial may add or remove at most
two vertices from $R$, at most $2 \times \MAXDEG$ vertices may be
added or removed from $T$. So it affects $|N(v) \cap T|$ by at
most $2 \times \MAXDEG$. Also every assignment to trials that
results $|N(v) \cap T| \ge k$ can be certified by the outcome of
$\MAXDEG k$ trials. By Talagrand's Inequality
$$\Prob[|N(v) \cap T| >\frac{\Delta}{1000}+ t+ 34800 \sqrt{\MAXDEG
\frac{\Delta}{1000}}]\le 4 e^{- 1000t^2 / 8 \times 4 \times
\MAXDEG^3 \Delta}.$$
Substituting $t = \frac{\Delta}{1000}$ in both inequalities
completes the proof of Claims 1 and 2.

\noindent {\bf Proof of Claim 3:} To prove this claim, instead of
directly applying Talagrand's Inequality to the random variable
$|N(v) \cap Q|$, we apply it to the random variable $X_v=|N(v)
\setminus Q|=\deg(v) -|N(v) \cap Q|$. Changing the outcome of each
trial affects $X_v$ by at most $2$, and every outcome of trials
that results $X_v \ge k$ can be certified by the outcome of
$\MINDEG k$ trials. Since $X_v=\deg(v) -|N(v) \cap Q|$, and so
$\Ex[X_v] \ge \deg(v) - \frac{\Delta}{1000}$, Talagrand's
Inequality implies that
$$\Prob[X_v < \deg(v)-\frac{\Delta}{1000}-t- 120 \sqrt{\MINDEG
\Delta}]\le 4e^{-t^2/8\times4\times20\times
(\deg(v)-\frac{\Delta}{1000})}\le 4e^{-t^2/1000\Delta}.$$
Substituting $t =\frac{\Delta}{1000}$ in the inequailty above
shows that
$$\Prob[X_v <\deg(v)-\frac{\Delta}{300}]\le \frac{1}{3\Delta^{7}}.$$
Since $X_v=\deg(v) -|N(v) \cap Q|$, we have
$$\Prob[|N(v) \cap Q| >\frac{\Delta}{300}]\le \frac{1}{3\Delta^{7}}.$$
\end{proof}

In Phase~{\bf II}, some of the edges that are incident to vertices
in $L$ will be uncolored. Since the other endpoints of these edges
may lie in $V(G')-V(H)-L$, to guarantee that these vertices will
remain distinguishable we need to prove a stronger condition than
just showing that they are distinguishable in $\sigma_1$. The
following lemma is needed to prove Lemma~\ref{phaseI}~(b) which
can imply that with positive probability all vertices in
$V(G')-V(H)-L$ will remain distinguishable after applying
Phase~{\bf II}.

\begin{lemma}
\label{triangle} For every two adjacent vertices $u,v \in
V(G')-V(H)$ where $\deg(u)=\deg(v)$,

$$\Prob[u \not\in L \wedge (|S_{\sigma_1}(u)
\bigtriangleup S_{\sigma_1}(v)| < 10)] < \frac{1}{\Delta^{7}}.$$
\end{lemma}
\begin{proof}
It is sufficient to prove that $\Prob[(|S_{\sigma_1}(u)
\bigtriangleup S_{\sigma_1}(v)| < 10) | u \not \in L] <
\frac{1}{\Delta^{7}}.$ Since we want to prove an upper bound for
$\Prob[(|S_{\sigma_1}(u) \bigtriangleup S_{\sigma_1}(v)| < 10) | u
\not \in L]$, we can assume that $S_\sigma(u)=S_\sigma(v)=S$.
(Remember that $\sigma$ is the edge coloring before applying
Phase~{\bf I}.)

Suppose that $\deg_{U_{\sigma_1}}(u) =k$ and $k \ge 20$. If
$|S_{\sigma_1}(u) \bigtriangleup S_{\sigma_1}(v)| < 10$, then
there are at least $k-10$ colors in $S-S_{\sigma_1}(u)$ which are
also in $S-S_{\sigma_1}(v)$. Since there are at most two edges
between $u$ and $v$, the probability of this occuring is at most
$${k \choose k-10}\left(\frac{\PrConst}{\Delta}\right)^{k-12}<\frac{1}{\Delta^{7}}.$$
So
$$\Prob[(|S_{\sigma_1}(u) \bigtriangleup S_{\sigma_1}(v)| < 10) | u
\not \in L] < \frac{1}{\Delta^{7}}.$$
\end{proof}

\begin{lemma}
\label{phaseI} If we apply Phase~{\bf I} to the edge coloring
$\sigma$, and obtain a partial edge coloring $\sigma_1$, then with
positive probability
\begin{enumerate}
\item[$(a)$:] For every vertex $v \in V(G')-V(H)$, we have $|N(v)
\cap L| \le \frac{\Delta}{\UDEG}$,

\item[$(b)$:] For every two adjacent vertices $u,v \in V(G')-V(H)$
such that $\deg(u)=\deg(v)$ and $u \not\in L$, we have
$|S_{\sigma_1}(u) \bigtriangleup S_{\sigma_1}(v)| \ge 10$.
\end{enumerate}
\end{lemma}
\begin{proof}
In order to prove the lemma, we define the following ``bad''
events
\begin{itemize}
\item For a vertex $v \in V(G')-V(H)$, let $A_v$ be the event that
$|N(v) \cap L| > \frac{\Delta}{\UDEG}$.

\item For every edge $uv$ such that $u,v \in V(G')-V(H)$,
$\deg(u)=\deg(v)$, let $A_{uv}$ be the event that  $u \not\in L$
and $|S_{\sigma_1}(u) \bigtriangleup S_{\sigma_1}(v)| < 10$.
\end{itemize}
It is easy to see that each event $A_X$ is mutually independent of
all events $A_Y$ such that all vertices in $X$ are in a distance
of at least $6$ from all vertices in $Y$. Hence each event is
independent of all events but at most $2\Delta^6$ events, and by
Lemmas~\ref{intersectL} and~\ref{triangle} each event occurs with
a probability of at most $\frac{1}{\Delta^{7}}$. So the Local
Lemma implies this lemma.
\end{proof}
\\
\\
\noindent {\bf Phase II}
\begin{enumerate} \item For every vertex
$u \in L$ choose $5$ edges $uv_i$ ($1 \le i \le 5$) such that $v_i
\in V(G')-L$, and $uv_i$ is not an unused edge uniformly at random
and uncolor them.
\end{enumerate}

Call this new partial edge coloring $\sigma_2$. For every vertex
$u \in V(G')-V(H)$, let $\UC'_u$ denote the set of the edges that
are incident to $u$ and are uncolored in Phase~{\bf II}.

\begin{lemma}
\label{phaseII} Suppose that $\sigma_1$ is a partial edge coloring
of $G'$ which satisfies Properties~$(a)$ and~$(b)$ in
Lemma~\ref{phaseI}. If we apply Phase~{\bf II} to $\sigma_1$, and
obtain a partial coloring $\sigma_2$, then with positive
probability we have
\begin{itemize}
\item[$(a)$:] For every vertex $v \in V(G')-V(H)-L$, $|\UC'_v| \le
4$.

\item[$(b)$:] For every two adjacent vertices $u,v \in
V(G')-V(H)$ such that $\deg(u)=\deg(v)$, we have $S_{\sigma_2}(u)
\neq S_{\sigma_2}(v)$.
\end{itemize}
\end{lemma}
\begin{proof}
Suppose that $u,v \in V(G')-V(H)$ are two adjacent vertices such
that $\deg(u)=\deg(v)$ and $u \not\in L$. Then since by
Property~$(b)$ in Lemma~\ref{phaseI} we know that
$|S_{\sigma_1}(u) \bigtriangleup S_{\sigma_1}(v)| \ge 10$, (a)
implies that $S_{\sigma_2}(u) \neq S_{\sigma_2}(v)$. It is trivial
that $|\UC'_u|=5$ for all vertices $u \in L$. So it is sufficient
to prove~(b$'$) instead of~(b), where

\begin{enumerate}
\item[$(b')$:]{\it For every two adjacent vertices $u,v \in L$
such that $\deg(u)=\deg(v)$, we have  $S_{\sigma_2}(u) \neq
 S_{\sigma_2}(v)$.}
\end{enumerate}

We will apply the Local Lemma to show that with positive
probability (a) and (b$'$) hold. We define the following two types
of ``bad'' events.
\begin{itemize}
\item For every vertex $u \in V(G')-V(H)-L$, and every five edges
$uv_1,uv_2,\ldots,uv_5$ such that $v_i \in L$ and $uv_i$ is not an
unused edge in $\sigma_1$, let $A_{u,\{v_1,\ldots,v_5\}}$ denote
the event that $uv_i \in \UC'_u$ for all $1 \le i \le 5$. Note
that since $G'$ is a multigraph, $\{v_1,\ldots,v_5\}$ is a
multiset and may have duplicate elements.

\item For every two adjacent vertices $u,v \in L$ where
$\deg(u)=\deg(v)$, let $A_{uv}$ denote the event that
$S_{\sigma_2}(u)=S_{\sigma_2}(v)$.
\end{itemize}

First we should give estimates of
$\Prob[A_{u,\{v_1,\ldots,v_5\}}]$ and $\Prob[A_{uv}]$. For every
vertex $v \in V(G')-V(H)$, we have $\deg(v) \ge \frac{\Delta}{3}$,
$|N(v) \cap L| \le \frac{\Delta}{\UDEG}$, and the unused degree of
$v$ is at most $\MAXDEG$. So every edge $uv$ is uncolored in
Phase~{\bf II} with a probability of at most
$\frac{5}{\left(\frac{\Delta}{3}-\frac{\Delta}{100}-290\right)} <
\frac{20}{\Delta}$, and also two parallel edges are uncolored with
a probability of at most $\frac{5}{\left(\frac{\Delta}{4}\right)}
\times \frac{4}{(\frac{\Delta}{4}-1)}$ which is less than
$(\frac{20}{\Delta})^2$. Since $u \not\in L$, this implies that

\begin{equation}
\Prob[A_{u,\{v_1,\ldots,v_5\}}] \le
\left(\frac{20}{\Delta}\right)^5.
\end{equation}

Suppose that $S_{\sigma_2}(u)=S_{\sigma_2}(v)$ for two adjacent
vertices $u,v \in L$ where $\deg(u)=\deg(v)$. Then since for each
of $u$ and $v$, five new incident edges are uncolored, we have
$\deg_{U_{\sigma_1}}(u)=\deg_{U_{\sigma_1}}(v)$. Since
$S_{\sigma_2}(u)=S_{\sigma_2}(v)$, all the five edges in $\UC'_v$
are determined by the five edges in $\UC'_u$, and since the edges
in $\UC'_v$ are chosen independently from the edges in $\UC'_u$,
we have

\begin{equation}
\Prob[A_{uv}]\le \frac{1}{{\frac{\Delta}{4} \choose 5}} \le
\left(\frac{20}{\Delta}\right)^5.
\end{equation}

Construct a graph $D$ whose vertices are all the events of the
above two types, in which two vertices $A_{u,\{v_1,\ldots,v_5\}}$
and $A_{u',\{v'_1,\ldots,v'_5\}}$ are adjacent if and only if
$\{v_1,\ldots,v_5\} \cap \{v'_1,\ldots,v'_5\} \neq \emptyset$
(even if $u=u'$), two vertices $A_{u,\{v_1,\ldots,v_5\}}$ and
$A_{vw}$ are adjacent if and only if $\{v_1,\ldots,v_5\} \cap
\{v,w\}\neq \emptyset$, and two vertices $A_{vw}$ and $A_{v'w'}$
are adjacent if $\{v,w\} \cap \{v',w'\} \neq \emptyset$. It is
easy to see that each vertex of $D$ is mutually independent of all
the vertices that are not adjacent to it. Since for every vertex
$w \in L$, there are at most $\Delta {\frac{\Delta}{100} \choose
4}$ events of the form $A_{u,\{w,v_2\ldots,v_5\}}$, and there are
at most $\frac{\Delta}{100}$ events of the form $A_{wv}$, the
maximum degree of $D$ is at most $$5 \left(\Delta
{\frac{\Delta}{100} \choose 4}+ \frac{\Delta}{100}\right) <
6\Delta {\frac{\Delta}{100} \choose 4}<\frac{6\Delta^5}{10^8}.$$
Then since $4\left(\frac{20}{\Delta}\right)^5
\frac{6\Delta^5}{10^8}<1$, the Local Lemma implies
Lemma~\ref{phaseII}.
\end{proof}

Let $\sigma_1$ and $\sigma_2$ be the partial colorings that are
guaranteed to exist by Lemmas~\ref{phaseI} and~\ref{phaseII}. The
maximum unused degree of $\sigma_1$ is $\MAXDEG$. So by
Lemma~\ref{phaseII}~(a) the maximum unused degree of $\sigma_2$ is
$294$. By Vizing's Theorem we can color the unused graph
$U_{\sigma_2}$ with $296$ new colors, and obtain a
$(\Delta+298)$-edge coloring $c$ of $G'$. By
Lemma~\ref{phaseII}~(b) for every two adjacent vertices $u,v \in
V(G')-V(H)$, we have $S_c(u) \neq S_c(v)$.

\subsection{Third Step}

In this step we begin with the edge coloring $c_1=c$, and
repeatedly modify it to eventually obtain a
$(\Delta+300)$-\avd-coloring, $c_{\final}$ of $G'$. We only
recolor the edges in $H$, and since the degree of every vertex in
$V(G')-V(H)$ is more than the degree of every vertex in $H$, all
of the vertices in $V(G')-V(H)$ will remain distinguishable.

Obtain $c_{k+1}$ from $c_k$ as in the following. Suppose that for
an edge $uv \in H$, $S_{c_k}(u)=S_{c_k}(v)$. Since $H$ does not
have any isolated edge,  we can assume that $\deg_H(u) \ge 2$. Let
$v_1,v_2, \ldots,v_r$ be the neighbors of $u$ in $H$, where
$r=\deg_H(u)$. Uncolor all edges $uv_1,uv_2,\ldots,uv_r$ to obtain
a partial edge coloring $c_k'$.

\begin{lemma}
\label{thirdstep} There exist sets $L(uv_i) \subseteq
\{1,\ldots,\Delta+300\}$ for $1 \le i \le r$, all of size $r+300$
such that
\begin{enumerate}
\item[$(a)$:] In the partial edge coloring $c_k'$ all the colors
in $L(uv_i)$ are available for the edge $uv_i$.

\item[$(b)$:] If $v' \neq u$ is adjacent to $v_i$ and
$S_{c_k'}(v') \setminus S_{c_k'}(v_i)$ has only one element $x$,
then $x \not\in L(uv_i)$.

\end{enumerate}
\end{lemma}
\begin{proof}
Since $\deg(u),\deg(v_i) \le \frac{\Delta}{3}$, there are at least
$\frac{\Delta}{3}+r+300$ available colors for $uv_i$. Furthermore
 there are at most $\frac{\Delta}{3}$ colors for $uv_i$ which
dissatisfy Condition (b). This implies that the desired sets
exist.
\end{proof}

Consider all {\em proper} completions of $c_k'$ such that the
color of each $uv_i$ is chosen from $L(uv_i)$. Note that the
number of these completions is at least $(r+300)(r+299)\ldots
301$. Pick one of them randomly and uniformly, and call it
$\hat{c_k}$. Since $H$ does not have any multiple edges, for every
$v_i$, there are at most $(r+300)(r-1)!$ different possible
choices of $\hat{c_k}$ for which $S_{\hat{c_k}}(u) =
S_{\hat{c_k}}(v_i)$. So we have
$$\Prob[\bigvee_{1\le i \le r} S_{\hat{c_k}}(u) =
S_{\hat{c_k}}(v_i)] \le \sum_{i=1}^{r} \Prob[S_{\hat{c_k}}(u) =
S_{\hat{c_k}}(v_i)] \le \sum_{i=1}^{r} \frac{(r-1)!}{(r+299)\ldots
301}<1.$$

This implies that there exists a completion $c_{k+1}$ of $c_k'$
such that $u$ is distinguishable, and by Lemma~\ref{thirdstep}~(b)
all neighbors of $u$ are also distinguishable. Since we only
changed the color of the edges incident to $u$, the number of
indistinguishable vertices is decreased. Hence by repeatedly
applying this procedure we will eventually obtain a
$(\Delta+300)$-\avd-coloring of $G'$.

\section*{Acknowledgements}
I want to thank my supervisor, Michael Molloy, for reading the
draft version of this paper and his valuable comments and
suggestions.

\bibliographystyle{plain}
\bibliography{avd}

\begin{thebibliography}{10}

\bibitem{Akbari}
S.~Akbari, H.~Bidkhori, and N.~Nosrati.
\newblock $r$-strong edge colorings of graphs.
\newblock {\em Discrete Math.}, to appear.

\bibitem{Balister2}
P.N. Balister, E.~Gy\"{o}ri, J.~Lehel, and R.H. Schelp.
\newblock Adjacent vertex distinguishing edge-colorings.
\newblock {\em J. Graph Theory}, to appear.

\bibitem{Balister}
P.N. Balister, O.M. Riordan, and R.H. Schelp.
\newblock Vertex-distinguishing edge colorings of graphs.
\newblock {\em J. Graph Theory}, 42:95--109, 2003.

\bibitem{Burris}
A.~C. Burris and R.~H. Schelp.
\newblock Vertex-distinguishing proper edge-colorings.
\newblock {\em J. Graph Theory}, 26 (2):73--82, 1997.

\bibitem{favaron}
O.~Favaron, H.~Li, and R.~H. Schelp.
\newblock Strong edge coloring of graphs,.
\newblock {\em Discrete Math.}, 159:103--109, 1996.

\bibitem{ghandehari}
M.~Ghandehari and H.~Hatami.
\newblock Two upper bounds for the strong edge chromatic number.
\newblock {\em Preprint}.

\bibitem{Hornak2}
M.~Hor\v{n}\'{a}k and R.~Sot\'{a}k.
\newblock Asymptotic behavior of the observability of {$Q_n$}.
\newblock {\em Discrete Math.}, 176:139--148, 1997.

\bibitem{Luczak}
M.~Karo\'{n}ski, T.~{\L}uczak, and A.~Thomason.
\newblock Edge weights and vertex colours.
\newblock {\em J. Combin. Theory Ser. B}, 91:151--157, 2004.

\bibitem{molloy}
M.~Molloy and B.~Reed.
\newblock {\em Graph Colouring and the Probabilistic Method}.
\newblock Springer, 2002.

\bibitem{west}
D.~B. West.
\newblock {\em Introduction to Graph Theory}.
\newblock Prentice-Hall, Inc, United States of American, 2001.
\newblock 2nd Edition.

\bibitem{Zhang}
Z.~Zhang, L.~Liu, and J.~Wang.
\newblock Adjacent strong edge coloring of graphs.
\newblock {\em Appl. Math. Letters}, 15:623--626, 2002.

\end{thebibliography}

\end{document}